\documentclass[10pt]{article}

\usepackage{amssymb,amsmath,latexsym}

\usepackage[latin1]{inputenc}
\usepackage{amsmath}
\usepackage{amsfonts}
\usepackage{amssymb}
\usepackage{makeidx}
\usepackage{bm}
\usepackage{color}
\usepackage{graphicx}
\usepackage{amscd}
\usepackage{fancyhdr}
\usepackage{latex8}

\newtheorem{theo.}{Theorem}[section]
\newtheorem{lem.}[theo.]{Lemma}
\newtheorem{cor.}[theo.]{Corollary}
\newtheorem{prop.}[theo.]{Proposition}
\newtheorem{def.}{Definition}
\newcommand{\fkfield}{\ensuremath{\mathbb{F}/\mathbb{K}}}

\newcommand{\kfield}{\ensuremath{\mathbb{K}}}
\newcommand{\oring}{\ensuremath{\mathcal{O}}}
\newcommand{\pideal}{\ensuremath{\mathcal{P}}}
\newcommand{\qideal}{\ensuremath{\mathcal{Q}}}
\newcommand{\ffield}{\ensuremath{\mathbb{F}}}
\newcommand{\zring}{\ensuremath{\mathbb{Z}}}

\newcommand{\zee}{\ensuremath{\mathbb Z}}

\newcommand{\places}{\ensuremath{\mathbb{P}_{\mathbb{F}}}}

\newcommand{\divisors}{\ensuremath{Div(\mathbb{\mathbb{F}})}}

\newcommand{\principal}{\ensuremath{\mathcal{P} _{F}}}

\newcommand{\jacobian}{\ensuremath{\mathbb{J}_{\mathbb{F}}}}

\begin{document}

\title{Explicit endomorphism of the Jacobian of a hyperelliptic function field of genus 2 using base field operations}

\author{Eduardo Ruiz Duarte\\
  Instituto de Matem\'{a}ticas (UNAM)\\
  \texttt{rduarte@ciencias.unam.mx}
\and
Octavio P\'{a}ez Osuna\\
  Passaic Co. Comm. College \\
  \texttt{opaezosuna@pccc.edu}
}

\date{\today}

\maketitle
\thispagestyle{empty}

\begin{abstract}

We present an efficient endomorphism for the Jacobian of a curve $C$ of genus 2 for divisors having a Non disjoint support. This extends the work of Costello in ~\cite{Costello} who calculated explicit formul\ae\space for divisor doubling and addition of divisors with disjoint support in $\jacobian(C)$ using only base field operations. Explicit formul\ae\space is presented for this third case and a different approach for divisor doubling. 

\end{abstract} 

\Section{Introduction}

High speed implementation of asymmetric cryptosystems is an important requirement to secure communications using devices with a reduced processor power, commonly RSA is used to provide asymmetric cryptography but this is not always the best solution for public key cryptography because of the high processing requirements for big random prime generation. \\
Elliptic curve cryptography (ECC) provides shorter keys and faster computation than RSA in some platforms and they offer the same security with less bits in the key ~\cite{faster}.

The Jacobian of a hyperelliptic curve is an excellent candidate for discrete logarithm problem based cryptography with the same safety benefits as RSA but with even shorter bit-lengths keys~\cite{hyperbits}.

Hyperelliptic curve cryptography (HECC) has been studied by Lange, Wollinger ~\cite{Lange}, ~\cite{Wollinger} to get explicit formul\ae\space for the calculation of the group operation over the Jacobian of genus 2 hyperelliptic curves, there are algorithms for arbitrary genus ~\cite{Cantor}, we will present a method to calculate the group operation using only arithmetic in the base field without polynomial pseudo inversion (Lange) using a geometrical approach for genus 2 hyperelliptic function fields using divisor theory and the Mumford representation of the elements of the Jacobian of the hyperelliptic function field, this approach has been studied by ~\cite{Costello} solving a system of linear equations over the base field to get the coefficients of the polynomial that defines the addition of two divisors in two cases: point doubling and regular addition, but regular addition has a subcase when the pair of divisors to be added share a point (place) in their supports, this third case is presented using formal differentiation with explicit formul\ae.

Following ~\cite{Stichtenoth}, an algebraic function field $\fkfield$ in one variable over $\kfield$ is an extension $\kfield\subseteq \ffield$ such that $\ffield$ is a finite algebraic extension of $\kfield(x)$ for some transcendental $x\in \ffield$ over $\kfield$.

A \textbf{valuation ring} in the function field $\fkfield$ is a ring $\oring\subseteq \ffield$ such that $\kfield \subsetneq \oring \subsetneq \ffield$ and for all $z\in \ffield$, $z\in\oring$ or $z^{-1}\in \oring$. These rings are \textbf{local rings} i.e. they have only one maximal ideal $\pideal=\oring \setminus \oring^{\times}$ , where $\oring^{\times}$ are the units of $\oring$ and all the ideals of $\oring$ are principal. It follows that every $0\neq z\in \ffield$ has a unique representation $z=t^{n}u$ for some $n\in\zring$ and $u\in\oring^{\times}$. A \textbf{place} $\pideal$ of $\fkfield$ is the maximal ideal of a valuation ring $\oring$ of $\fkfield$. 

If $\oring$ is a valuation ring of $\fkfield$ and $\pideal$ is its maximal ideal, then $\oring$ is determined only by $\pideal$ and we denote this ring by: \\
 \\
$\oring_{\pideal} := \lbrace z\in\ffield \mid z^{-1}\notin \pideal \rbrace$ \\
 \\
Here we say that $\oring_{\pideal}:=\oring$ is the valuation ring at the place $\pideal$ and the number of elements of $\places$ is infinite ~\cite{Stichtenoth}.

Let $\pideal \in \places$, the map $v_{\pideal}:\ffield \mapsto \zring\cup\lbrace \infty \rbrace$ is a discrete valuation of $\fkfield$ associated to $\pideal$ in the way that if $t$ is a \textbf{uniformization variable} then for all $0\neq z \in \ffield$ exists a unique representation of $z$, $z=t^{n}u$ with $u\in\oring_{\pideal}^{\times}$ and $n\in\zring$, then we define $v_{\pideal}(z):=n$ and $v_{\pideal}(0):=\infty$

$\mathbb{F}_{\pideal}:=\mathcal{O}_{\pideal}/\pideal$ will be the field of residual classes of $\pideal$, (this is field because $\pideal$ is maximal in $\mathcal{O}_{\pideal}$, this classes will be defined as $x+\pideal:=x(\pideal)$.
$deg(\pideal):=[\mathbb{F}_{\pideal}:\mathbb{K}]$ will be the degree of $P$. It follows that $deg(\pideal)$ is finite.

Let $z\in\ffield$ and $\pideal \in \places$, we say that $\pideal$ is a zero of $z$ if $v_{\pideal}(z)>0$ and $\pideal$ is a pole of $z$ if $v_{\pideal}(z)<0$ this is that if $v_{\pideal}(z)=m>0$ then $\pideal$ is a zero of $z$ of order m, if $v_{\pideal}(z)=-m<0$ we say that $\pideal$ is a pole of $z$ of order $m$. Every $z\in\ffield$ transcendental over $\kfield$ has at least one zero and one pole, in fact it has the same finite number of zeroes and poles.

A divisor is a finite formal sum of places
\begin{displaymath}
D=\sum_{\pideal\in\mathbb{\pideal}_{\mathbb{F}}}{n_{\pideal} \pideal} \mbox{ with } n_{\pideal}\in\zee ,\mbox{ and almost all }  n_{P}=0.
\end{displaymath}
The support of a divisor $D\in\divisors$ is defined as
\begin{displaymath}
supp(D):=\{ \pideal\in\mathbb{\pideal}_{\mathbb{F}} | n_{\pideal}\neq 0\} .
\end{displaymath}
Given $D=\sum{n_{\pideal} \pideal}$ y $D'=\sum{n_{\pideal}' \pideal}$ the sum is done coefficient wise:
        \begin{displaymath}
        D+D'=\sum_{\pideal\in\mathbb{\pideal}_{\mathbb{F}}}{(n_{\pideal}+n_{\pideal}') \pideal} .
        \end{displaymath}

        The zero element of $\divisors$ is:
        \begin{displaymath}
0:=\sum_{\pideal\in\mathbb{\pideal}_{\mathbb{F}}}{n_{\pideal} \pideal} \mbox{ with all the } n_{P}=0.
\end{displaymath}

For $\qideal\in\mathbb{P}_{\mathbb{F}}$ and $D\in\divisors$ we define $v_{\qideal}(D)=n_{\qideal}$,
    then
     \begin{displaymath}
     supp(D)=\{ \pideal\in\mathbb{P}_{\mathbb{F}} | v_{\pideal}(D)\neq 0\} \mbox{ and }
     D=\sum{v_{\pideal}(D)\cdot
             \pideal} .
             \end{displaymath}

A partial order in the group of divisors is given by
\begin{displaymath}
D_{1}\leq D_{2}:\Longleftrightarrow v_{\pideal}(D_{1})\leq v_{\pideal}(D_{2}) \mbox{
        for all } \pideal\in\mathbb{\pideal}_{\mathbb{F}} .
        \end{displaymath}

        A divisor such that $D\geq 0$ is \textbf{positive} or \textbf{effective} .
        The \textbf{degree of a divisor} is defined as
        \begin{displaymath}
        \partial(D):=\sum_{\pideal\in\mathbb{P}_{\mathbb{F}}}{v_{\pideal}(D)\cdot deg(\pideal)}
        \end{displaymath}

For a function $0\neq x\in \ffield$ let $Z$ be the set of zeros and $N$ be the set of poles of
        $x$ in $\mathbb{P}_{\mathbb{F}}$. we define
        \begin{displaymath}
        (x)_{0}:=\sum_{\pideal\in Z}{v_{\pideal}(x)\pideal} , \mbox{ zero divisor of } x, 
        \end{displaymath}
        \begin{displaymath}
        (x)_{\infty}:=-\sum_{\pideal\in N}{v_{\pideal}(x)\pideal} , \mbox{ pole divisor of } x, 
        \end{displaymath}
        \begin{displaymath}
        (x):=(x)_{0}-(x)_{\infty} ,\mbox{ principal divisor of} x. 
        \end{displaymath}
Divisors $(x)_{0}$ and $(x)_{\infty}$ are effective divisors, and
        \begin{equation}
        (x)=\sum_{\pideal\in\mathbb{P}_{\mathbb{F}}}{v_{\pideal}(x)\pideal}
        \end{equation}

The set
\begin{displaymath}
\principal :=\{ (x)|0\neq x \in F\}
\end{displaymath}
is called the \textbf{subgroup of principal divisors} of $\mathbb{F}/\mathbb{K}$.

The quotient group
\begin{displaymath}
\jacobian := \divisors^0 / \principal
                  \end{displaymath}
                  will be defined as \textbf{the group of divisor classes or Jacobian of $\fkfield$}. For $D\in \divisors$, the corresponding element in $\jacobian$ is denoted by $[D]$, the
                  class of $D$. Two divisors $D$, $D'\in \divisors$ are equivalent
                  ($D\sim D'$) if $[D]=[D']$, this is that, $D=D'+(x)$ for some $x\in
                  \mathbb{F}\backslash \{ 0\}$. This is an equivalence relation.

\section{Hyperelliptic function fields}

A hyperelliptic function field over $\kfield$ is a function field $\fkfield$ with genus $g\geq 2$ that contains a rational subfield $\kfield(x)\subseteq \ffield$ with $[\ffield:\kfield(x)]=2$

\begin{lem.}
\begin{enumerate}
\item{A function field $\fkfield$ of genus $g\geq 2$ is hyperelliptic if and only if there is a divisor $A\in\divisors$ with $\partial(A)=2$ and the dimension of the Riemann space at $A$ is greater or equal than 2}
\item{Every $\fkfield$ of genus 2 is hyperelliptic}
\end{enumerate}
\end{lem.}

\begin{theo.}
Let $D\in\divisors$ with $\partial(D)=0$, then there is a divisor $D'-rP\in[D]$ with $D'\geq 0$, $\partial(D')=r\leq g$ and $P$ a place. 
\end{theo.}
We will call the divisor of the previous Theorem the reduced divisor of $[D]$.

\begin{cor.}
Every element $[D]\in\jacobian$ of genus 2 with $D\equiv (p_x,p_y)+(q_x,q_y)-2\infty$, this divisor can be represented by the pair of functions $<u(x),v(x)>$ such that $u(p_x)=u(q_x)=0$, $v(p_{x})=p_{y}$ and $v(q_x)=q_y$ with $u$ monic, $deg(u)=g=2$ and $deg(v)=g-1=1$. The pair $<u(x),v(x)>$ is called \textbf{Mumford representation} of $[D]$. 
\end{cor.}

The next theorem will justify the closure of our method to do arithmetic in the Jacobian of a hyperelliptic curve. \\
\begin{theo.}
\textbf{Artin's aproximation theorem} ~\cite{Stichtenoth} \\ Let $\fkfield$ a function field and $\pideal_{1},\pideal_{2},...,\pideal_{n}\in\places$ different places in pairs of $\fkfield$, $x_{1},x_{2},...,x_{n}\in\ffield$ and $r_{1},r_{2},...,r_{n}\in\zring$ then there exists $x\in\ffield$ such that: \\
 \\
$v_{\pideal_{i}}(x-x_{i})=r_{i}$ with $i=1,2,...,n$
\end{theo.}
This theorem generalize the chinese remainder theorem and it will assure the existence of a curve that passes through the given points (places) with any degree of multiplicity at the hyperelliptic curve (valuation at the point)

\Section{Explicit addition on $\jacobian$ with Mumford divisors over a hyperelliptic function field of genus 2}
\SubSection{Case 1: $[D_{1}]\oplus [D_{2}]$ with $Supp(D_{1})\cap Supp(D_{2})=\emptyset$}

We just justify the structure of a hyperelliptic function field of genus two, we have the function field $\mathbb{K}(x,y)$ is such that $y^2=f(x)$ with $deg(f(x))=5$. \\
So we define the hyperelliptic curve as $C(x,y)=y^2-f(x)$ \\
 \\ 
Given two divisors $D_{1}=\pideal_{1}+\pideal_{2}-2\qideal_{\infty}$ and $D_{2}=\pideal_{1}'+\pideal_{2}'-2\qideal_{\infty}$ , we want to find the divisor class $[\pideal_{1}+\pideal_{2}-2\qideal_{\infty}]\oplus [\pideal_{1}'+\pideal_{2}'-2\qideal_{\infty}]$, to find this, we can use the approximation theorem to be sure of its existence, we have that there is a function  $L\in \mathbb{K}(x,y)$ and its principal divisor $(L)$ that has in the support the places of degree one $\pideal_{1},\pideal_{2},\pideal_{1}',\pideal_{2}',\bar{\pideal_{1}},\bar{\pideal_{2}}$, for this, we have to find an interpolation polynomial that passes through this places, so we make this polynomial equal to $y^2=f(x)$ $(deg(f(x))=5)$  and we solve for $\bar{\pideal_{1}},\bar{\pideal_{2}}$ and finally we make hyperelliptic involution to find $-[\bar{\pideal_{1}}+\bar{\pideal_{2}}-2\qideal_{\infty}]=[\pideal_{1}''+\pideal_{2}''-2\qideal_{\infty}]$. \\

The geometric intuition of what we want to find (blue) given two divisors with disjoint supports: \\
\begin{center}
\includegraphics[scale=0.5]{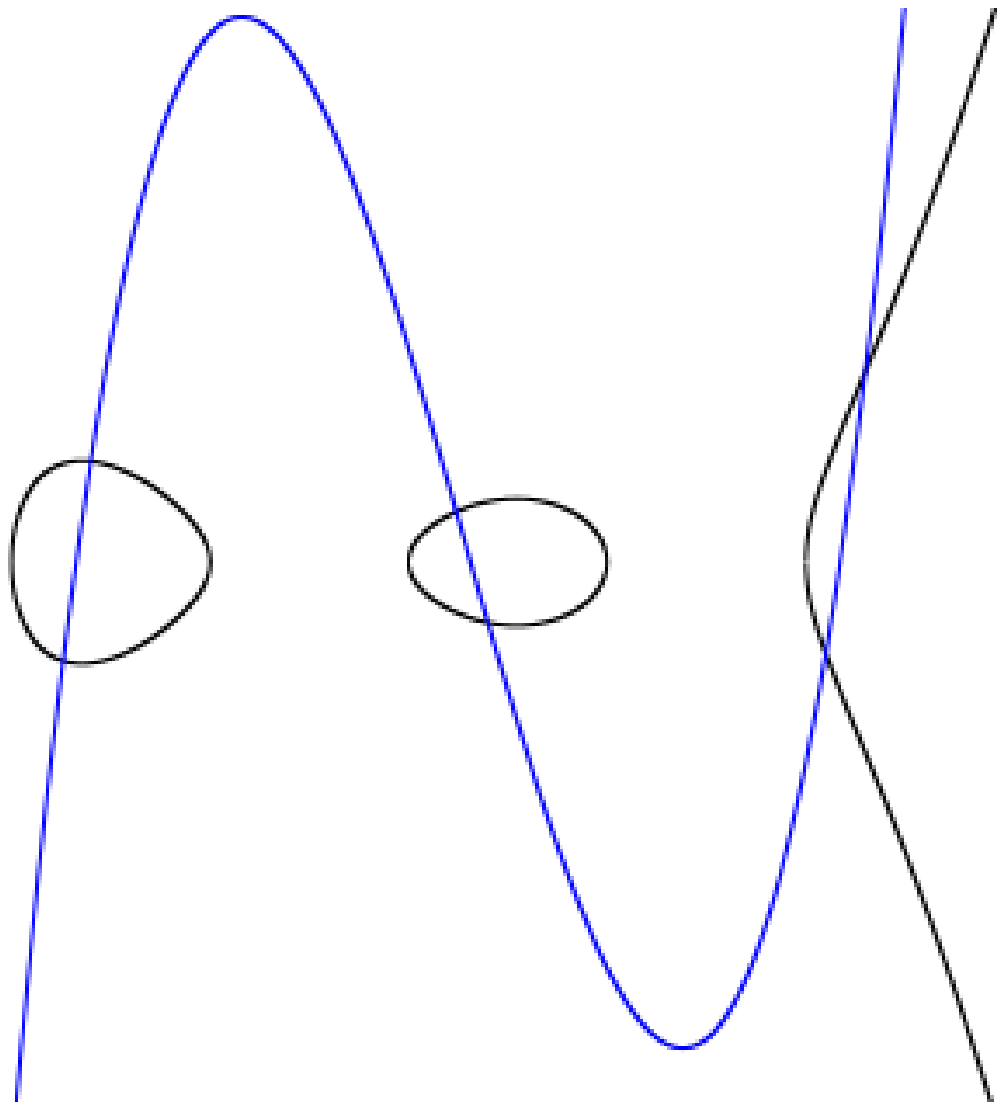}
\end{center}

Using Mumford representation, if $D$ is a reduced divisor as in the corresponding theorem, and we denote the places as $\pideal=(x,y)\in Supp(D)$, then  $D=(u,v)$ with $u(x)=0$ and $v(x)=y$ for all $\pideal=(x,y)\in Supp(D)$ so $deg(u)=2$ and $deg(v)=1$, so the addition is done as following: \\
 \\
$D_{1}=(u=x^2+ax+b,v=cx+d)$\\
$D_{2}=(u'=x^2+Ax+B,v'=Cx+D)$\\
 \\
We want to find $D_{3}=(u''=x^2+\alpha x+\beta,v''=\gamma x+\delta)$ as in the previous figure to represent $P_{1}'',P_{2}''\in Supp(D_{3})$, this results inverting $\bar{P_{1}}$ and $\bar{P_{2}}\in Supp(L)$ \\
 \\
To find the other elements of $Supp(L)$ $\bar{P_{1}}$ and $\bar{P_{2}}$ we have to find the interpolation polynomial for the given places, square it and then making it equal to $y^2=f(x)$: \\
 \\
We have that: \\
$L(x)=px^3+qx^2+rx+s$ \\
 \\

If we solve: \\
 \\

$L(x)-v(x)\equiv 0$ mod $u(x)$\\ 
 
$L(x)-v'(x)\equiv 0$ mod $u'(x)$\\
 \\
With $deg(u)=deg(u')=2$ we will have:\\
 \\
$R_{1}x+R_{2}$ \\
$R_{3}x+R_{4}$ \\
 \\
So $R_{i}=0$ and we will have a $4\times 4$ system of equations, the solutions are going to be the coefficients $p,q,r,s$ of $L(x)$, if we do the calculations reducing $L(x)-v(x)$ modulo $u(x)$ and $L(x)-v(x)$ modulo $u'(x)$ we can find that $r_{i}=0$ in general induces a matrix \\
 \\
$(px^3+qx^2+rx+s)-(cx+d) \equiv x(p(a^2-b)-qa+r-c)+p(ab)-qb+s-d$ mod $x^2+ax+b$ \\ 
 \\
$(px^3+qx^2+rx+s)-(Cx+D) \equiv x(p(A^2-B)-qA+R-C)+p(AB)-qB+s-D$ mod $x^2+Ax+B$  \\
 \\
This induces 4 equations: \\
 \\
$R_{1} = p(a^2-b)-qa+r-c$ \\
$R_{2} = p(ab)-qb+s-d$ \\
$R_{3} = p(A^2-B)-qA+R-C$ \\
$R_{4} = p(AB)-qB+s-D$ \\
 \\
As we know the values $a,b,c,d,A,B,C,D$ and we want $R_{i}=0 \quad \forall 1\leq i \leq 4$ to find the coefficients of $L(x)$ the system is this: \\
\[
\left[\begin{array}{cccc|c}
   a^2-b	&	-a	&	1	&	0	&	c\\
   ab		&	-b	&	0	&	1	&	d\\
   A^2-B	&	-A	&	1	&	0	&	C\\
   AB		&	-B	&	0	&	1	&	D 
\end{array}\right].
\]
The solution of this matrix give us the coefficients $p,q,r,s$ of $L(x)$ having this we just have to make it equal to the hyperelliptic curve $y^2=f(x)$: \\

$\frac{L(x)^2-f(x)}{u(x)u'(x)}=u''(x)$ \\

This happens because $deg(L)=6$ and it has $u$ and $u'$ as factors and the coordinates of $x$ are roots of the places of $Supp(D_{1})$ and $Supp(D_{2})$ so $u''(x)$ is the polynomial of degree 2 resulting of the quotient and this is going to have as roots the coordinates of $x$ in $Supp(D_{3})$. \\

To find $v''(x)$ we have that $deg(v'')=1$ we just have to evaluate the roots of $u''(x)$ in the hyperelliptic curve $C(x,y)=y^2 - f(x)$ over $\mathbb{K}$ but this can be a problem (computing roots) so another way is to check that: \\
 \\
$L(x) \equiv v''(x) \bmod u''(x)$
 \\
With this we have $[D_{1}]\oplus [D_{2}]= [D_{3}]=(u''(x),v''(x))$.  

\SubSection{Case 2: 2[D]}

See  ~\cite{Costello} but here we will show a slightly different approach, which in this case is a particular case of the previous, suppose we have $D=\mu+\omega-2\infty$ and we want to calculate $2[D]$, this divisor is linearly equivalent to a divisor
with \textbf{prime places of degree two}: $2[D] \sim [2\mu -2\infty] \oplus [2\omega - 2\infty]$ what we have here is $v_{\mu}(L)=2$ and $v_{\omega}(L)=2$ , the geometric intuition is that
we have two points which are 'tangent' to the hyperelliptic curve with degree of intersection two in both: \\
 \begin{center}
\includegraphics[scale=0.5]{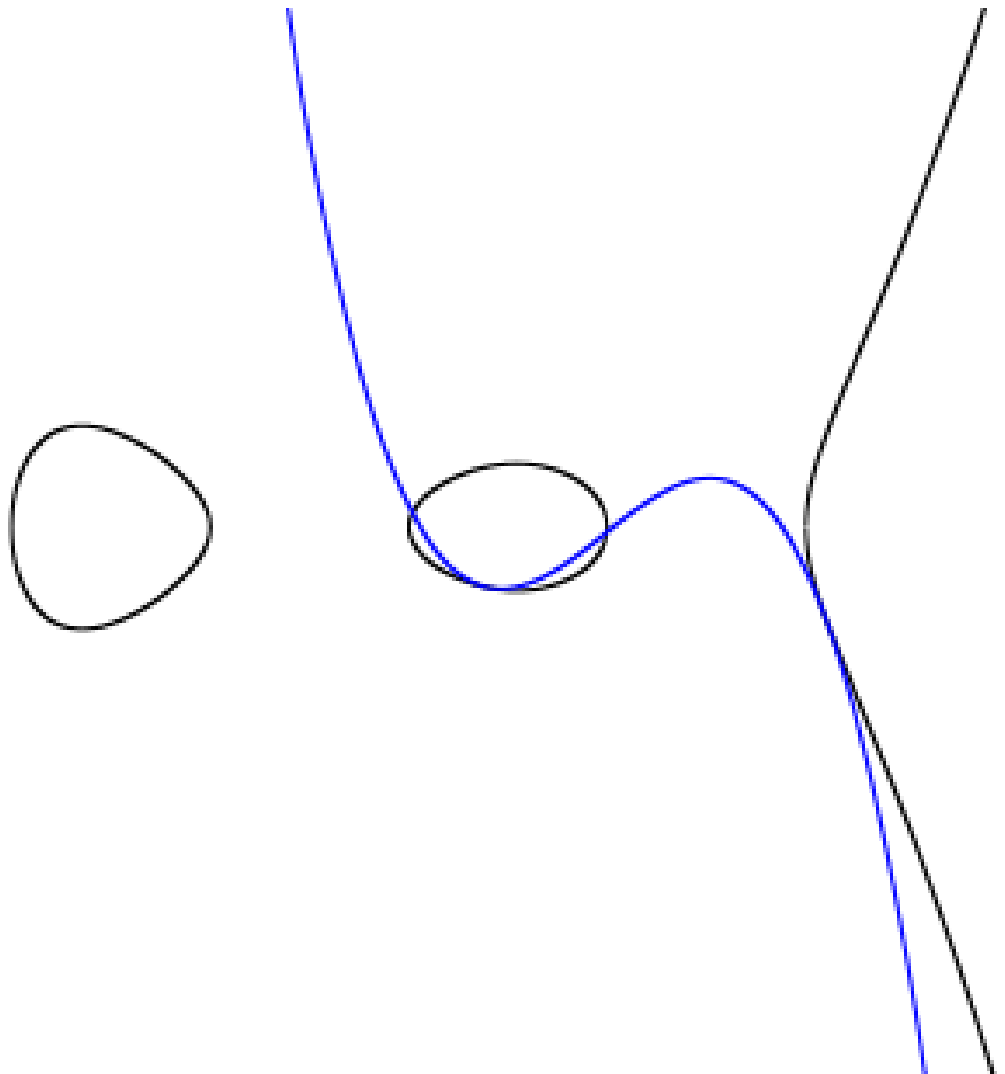}
\end{center}

The Mumford divisors are the following: \\
 \\

Let $\mu=(\mu_x,\mu_y)$ and $\omega=(\omega_x,\omega_y)$ be two points (places) of the hyperelliptic curve $C$ of genus two, $D=\mu+\omega-2\infty$ , we will calculate \\
 \\
 $2[D] = [D_1]\oplus [D_2] = [2[\mu-2\infty]\oplus [2\omega-2\infty]$ \\
 \\
$[D_1]=<x^2-2x\mu_x+{\mu_x}^2,\frac{dC}{dx}(\mu_x,\mu_y)x-\frac{dC}{dx}(\mu_x,\mu_y)\mu_x+\mu_y>$ \\
$[D_2]=<x^2-2x\omega_x+{\omega_x}^2,\frac{dC}{dx}(\omega_x,\omega_y)x-\frac{dC}{dx}(\omega_x,\omega_y)\omega_x+\omega_y>$ \\
 \\
Here we built the quadratic polynomial as a double root in the $x$ coordinate of the prime place of each divisor, and we used formal differentiation 
for the linear part to get the 'tangent' line to the curve $C$ at the given point (place). \\
 \\
So, using the matrix in the case 1 we need to solve for $P,Q,R,S$ such that $L(x)=Px^3+Qx^2+Rx+S$: 
 \\
$a=-2\mu_x$ \\
$b={\mu_x}^2$ \\
$c=\frac{dC}{dx}(\mu_{x},\mu_{y})$ \\
$d=\mu_{y}-\frac{dC}{dx}(\mu_{x},\mu_{y})\mu_{x}$ \\
 \\
$A=-2\omega_x$ \\
$B={\omega_x}^2$ \\
$C=\frac{dC}{dx}(\omega_{x},\omega_{y})$ \\
$D=\omega_{y}-\frac{dC}{dx}(\omega_{x},\omega_{y})\omega_{x}$

\[
\left[\begin{array}{cccc|c}
3{\mu_x}^2 & 2\mu_x & 1 & 0 & \frac{dC}{dx}(\mu_x,\mu_y) \\
-2{\mu_x}^{3} & -{\mu_x}^2 & 0 & 1 & \mu_{y}-\frac{dC}{dx}(\mu_x,\mu_y)\mu_x \\  
3{\omega_x}^2 & 2\omega_x & 1 & 0 & \frac{dC}{dx}(\omega_x,\omega_y) \\
-2{\omega_x}^{3} & -{\omega_x}^2 & 0 & 1 & \omega_{y}-\frac{dC}{dx}(\omega_x,\omega_y)\omega_x
\end{array}\right].
\]
In the same way with the solution of this system with get the coefficients of $L(x)$ and the new places to do hyperelliptic involution, so we are ready to compute $2[D]=(u'(x),v'(x))$ \\

\SubSection{Case 3: $[D_{1}]\oplus [D_{2}]$ with $Supp(D_{1})\cap Supp(D_{2})\neq\emptyset$}

in this new case we have that both divisors share a place $P=(s,t)$, let $D_1=P+\mu-2\infty$ and $D_2=P+\omega-2\infty$, We will need to detect explicitly the repeated place $P$ to calculate this case given the divisors in the Mumford 
notation, but this is easy because if $[D_{1}]=<u_1(x)=x^2+ax+b,v_1(x)>$ and $[D_{2}]=<u_2(x)=x^2+\alpha x+\beta,v_2(x)>$ the matrix generated by case 1 must be singular, so $u_1(x)$ and $u_2(x)$ have a common root so $u_1(x)=u_2(x)$ implies that the repeated root $x$ is $s=\frac{\beta-b}{a-\alpha}$, having this we have that $t=v_1(s)$ and the other places $\mu$ and $\omega$ can be calculated directly, so we are ready to calculate the addition in this case.

As in the case 2, we have that: \\
 \\
$[D_1]\oplus [D_2] \sim [2P -2\infty] \oplus [\mu+\omega-2\infty]$  \\
 \\
Geometrically we can sketch this situation as:
\begin{center}
\includegraphics[scale=0.5]{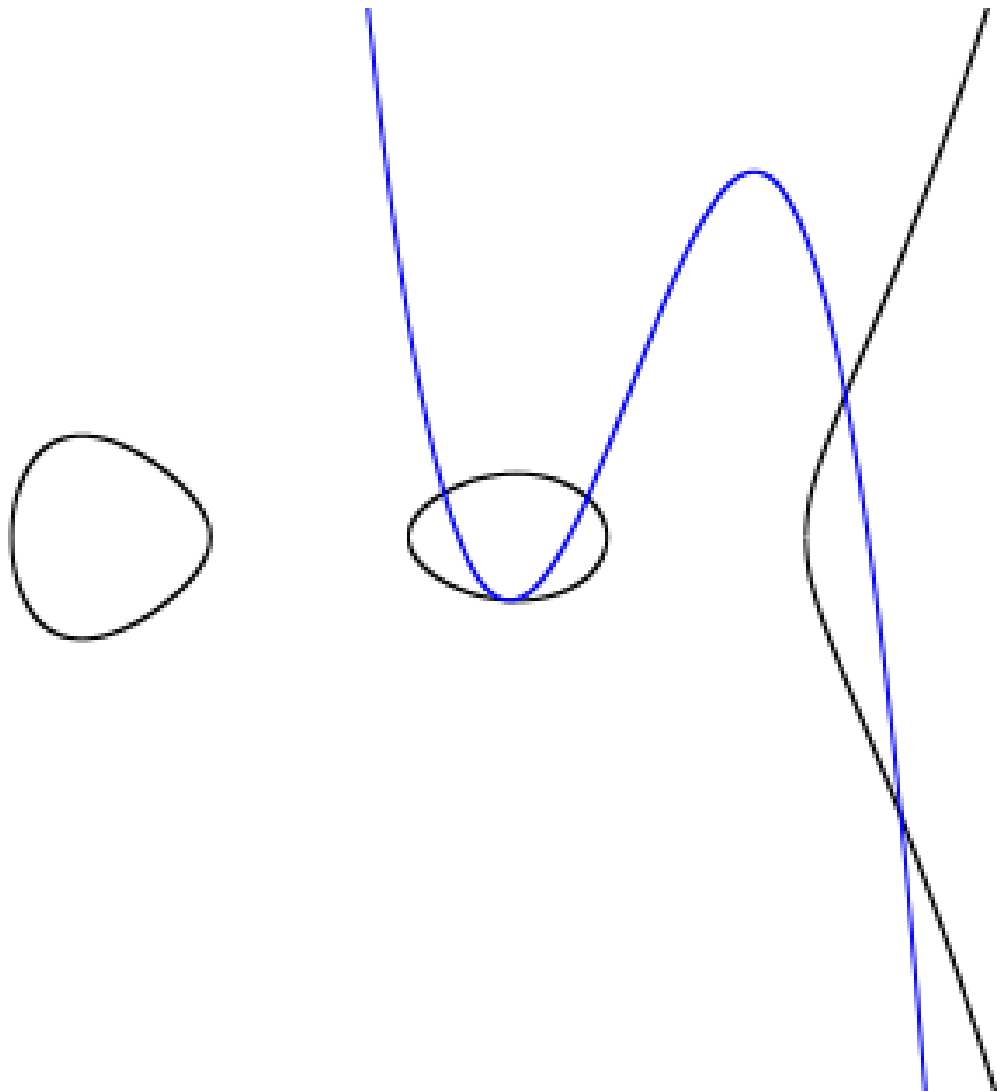}
\end{center}

We see here the place $P$ of degree 2 (blue curve tangent point at $C$). \\
 \\
So we rewrite $D_1$ and $D_2$ as:
 \\
$D_1=2(s,t)-2\infty$ \\
$D_2=(\mu_x,\mu_y)+(\omega_x,\omega_y)-2\infty$ \\
 \\
And we have the Mumford divisors in the following way:  \\
 \\
$[D_1]=<x^2-2sx+s^2,\frac{dC}{dx}(s,t)x-\frac{dC}{dx}(s,t)s+t>$ \\
$[D_2]=<x^2-x(\mu_x+\omega_x)+\mu_x\omega_x,\frac{\mu_{y}-\omega_{y}}{\mu_{x}-\omega_{x}}x+\frac{\mu_{x}\omega_{y}-\omega_{x}\mu_{y}}{\mu_{x}-\omega_{y}}>$ \\
 \\
Here we have that $[D_1]$ has a multiple root at $s$ and the linear part is the tangent line at $(s,t)$ over $C$ and $[D_2]$ the linear part is the line through $(\mu_x,\mu_y)$ and $(\omega_x,\omega_y)$, and the quadratic part has $\mu_x$ and $\omega_x$ as roots. \\
We use the case 1 and we have that: \\
As in the other cases we use the matrix of case 1 to obtain the coefficients of the interpolation polynomial $L\in \mathbb{K}(x,y)$ with $y^2=f(x)$ , $deg(f(x))=5$ and $L(x)=Px^3+Qx^2+Rx+S$.
so we solve for $P,Q,R,S$ the following matrix: \\

\[
\left[\begin{array}{cccc|c}
3s^2 & 2s & 1 & 0 & \frac{dC}{dx}(s,t) \\
-2s^3 & -s^2 & 0 & 1 & t-\frac{dC}{dx}(s,t)s \\
{\mu_{x}}^2+\omega_{x}(\mu_{x}+\omega_{x}) & \mu_{x}+\omega_{x}& 1 & 0 & \frac{\mu_{y}-\omega_{y}}{\mu_{x}-\omega_{x}} \\
-{\mu_{x}}^{2} \omega_{x}-\mu_{x}{\omega_{x}}^{2} & -\mu_{x}\omega_{x}&0&1&\frac{\mu_{x}\omega_{y}-\omega_{x}\mu_{y}}{\mu_{x}-\omega_{y}}
\end{array}\right].
\]


\SubSection{Conclusions}

Jacobians of hyperelliptic curves of genus 2 are a good candidate for asymmetric cryptography, so optimization of its endomorphism is an important task, the work of Costello ~\cite{Costello} is very important because it shows that the calculation of the addition in $<\jacobian,\oplus>$ can be done solving a system of linear equations over the base field in two cases given Mumford divisors of a genus 2 curve Jacobian, we have extended with a similar approach a third case when both divisors to add in the Jacobian share a point, the existence of the solution is backed up by the Artin's approximation theorem, the figures shown in this work were calculated using these formul\ae \space.


\end{document}